\documentclass[12pt]{article}

\usepackage{diagbox}
\usepackage{mathtools}
\usepackage{bbm}
\usepackage{latexsym}
\usepackage{epsfig}
\usepackage{amsmath,amsthm,amssymb,enumerate,hyperref}
\usepackage[active]{srcltx}
\usepackage{color}

\newcounter{rot}

\def\e{\epsilon}    
  \def\K{\Kappa}
     
   \def\p{\pi}
   
 \def\om{\omega}

\newtheorem{theorem}{Theorem}

\newcommand{\set}[1]{\left\{#1\right\}}

\def\Pr{\mbox{{\bf Pr}}}

\newcommand{\ignore}[1]{}

\newcommand{\C}[2]{C^{(#1)}_{#2}}

\usepackage{tikz}
\usetikzlibrary{arrows}
\usetikzlibrary{decorations}
\usetikzlibrary{shapes.misc}
\newcommand{\pic}[1]{
\begin{tikzpicture}
#1
\end{tikzpicture}}

\begin{document}
\author{Alan Frieze\thanks{Research supported in part by NSF grant DMS}\\Department of Mathematical Sciences\\Carnegie Mellon University\\Pittsburgh PA 15213}

\title{A note on spanning $K_r$-cycles in random graphs}
\maketitle
\begin{abstract}
We find a threshold for the existence of a collection of edge disjoint copies of $K_r$ that form a cyclic structure and span all vertices of $G_{n,p}$. We use a recent result of Riordan to give a two line proof of the main result.
\end{abstract}
\section{Introduction}
In a seminal paper, Johansson, Kahn and Vu \cite{JKV} solved the long standing open question of determining the threshold for the existence of $H$-factors in random graphs and hypergraphs. For some questions, the proof for hypergraphs turns out to be somewhat simpler than that of the related question in graphs. More precisely, the proof of the existence of a perfect matching in a random $r$-uniform hypergraph is simpler than the proof of the existence of a $K_r$-factor in $G_{n,p}$. Recently Riordan \cite{Rio} showed that one can avoid the more complicated proofs. He does this by proving a coupling between graphs and hypergraphs that enables one to infer graph factor thresholds from hypergraph matching thresholds. The aim of this short note is to show how to use this coupling to prove thresholds for some other spanning subgraphs. 

We are given a graph $G$ with $n$ vertices and an integer $r\geq 3$ where $n=(r-1)m$, $m$ integer. A {\em $K_r$-cycle} is a sequence $H_1,H_2,\ldots,H_m$ of copies of $K_r$ where (i) $V(H_i)\cap V(H_{i+1})=\set{v_i},i=1,2 \ldots,m$ ($v_{m+1}=v_1$ here) and (ii) $H_i$ and $H_j$ are vertex disjoint for $i\neq j$.
\begin{center}
\def\C4{
\draw [fill=black] (0,0) circle [radius=.05];
\draw [fill=black] (0,1) circle [radius=.05];
\draw [fill=black] (0.5,0.5) circle [radius=.05];
\draw [fill=black] (-0.5,0.5) circle [radius=.05];
\draw (0,0) -- (0.5,0.5);
\draw (0,0) -- (-0.5,0.5);
\draw (0,1) -- (0.5,0.5);
\draw (0,1) -- (-0.5,0.5);
\draw (0,0) -- (0,1);
\draw (-0.5,0.5) -- (0.5,0.5);
}
\pic{
\begin{scope}[shift={(0,0)}, scale=1]
\C4
\end{scope}
\begin{scope}[shift={(1,0)}, scale=1]
\C4
\end{scope}
\begin{scope}[shift={(2,0)}, scale=1]
\C4
\end{scope}
\begin{scope}[shift={(3,0)}, scale=1]
\C4
\end{scope}
\begin{scope}[shift={(3,1)}, scale=1]
\C4
\end{scope}
\begin{scope}[shift={(3,2)}, scale=1]
\C4
\end{scope}
\begin{scope}[shift={(2,2)}, scale=1]
\C4
\end{scope}
\begin{scope}[shift={(1,2)}, scale=1]
\C4
\end{scope}
\begin{scope}[shift={(0,2)}, scale=1]
\C4
\end{scope}
\begin{scope}[shift={(0,1)}, scale=1]
\C4
\end{scope}
\node at (1.5,-1) {A $K_4$-cycle};
}
\end{center}

We will prove the following theorem:
\begin{theorem}\label{th1}
$p=n^{-2/r}\log^{1/\binom{r}{2}}n$ is a threshold for $G_{n,p}$ to contain a spanning $K_r$-cycle.
\end{theorem}
\section{Proof of Theorem \ref{th1}}
For the proof, we need two results: the first will be Theorem 1 of Riordan \cite{Rio} combined with Theorem 2 of Heckel \cite{H}. 
\begin{theorem}\label{thrio}
Let $r\geq 3$ be given. There is a positive constant $\e$ such that if $p\leq n^{-2/r+\e}$ then, for some $\p\approx p^{\binom{r}{2}}$, we may couple $G=G_{n,p}$ and the random $r$-uniform hypergraph $H=H_{n,\p;r}$ such that w.h.p. to every edge $e$ of $H$ there is a corresponding copy of $K_r$ in $G$ with $V(K_r)=e$.
\end{theorem}
We will also need the following theorem from Dudek, Frieze, Loh and Speiss \cite{DFLS}, which removed some divisibility constraints from \cite{DF1}, \cite{F1}. A {\em loose} Hamilton cycle $C$ in an $r$-uniform hypergraph $H=(V,\mathcal{E})$ of order $n$ is a collection of edges of $H$ such that for some cyclic ordering of $V$, every edge consists of $r$ consecutive vertices, and for every pair of consecutive edges $E_{i-1},E_i$ in $C$ (in the natural ordering of the edges), we have $|E_{i-1}\cap E_i|=1$. 
\begin{theorem}\label{Ham}
Suppose $k\geq 3$.  If $\p=\om n^{1-r} \log n$ for $\om=\om(n)\to\infty$, then
\[
\lim_{\substack{n\to \infty\\(r-1) | n}}\Pr\left(H_{n,\p;r}\ contains\ a\ loose\ Hamilton\ cycle\right)=1.
\]
\end{theorem}
{\bf Proof of Theorem \ref{th1}}\\
First suppose that $p=\om n^{-2/r}\log^{1/\binom{r}{2}}n$. We couple $G_{n,p}$ with the hypergraph $H_{n,\p;r}$ as promised by Theorem \ref{thrio}. Because $p^{\binom{r}{2}}=(\om n^{-2/r}\log^{1/\binom{r}{2}}n)^{\binom{r}{2}}=\om^{\binom{r}{2}}n^{1-r}\log n$ we see from Theorem \ref{Ham}  that w.h.p. $H_{n,\p;r}$ contains a loose Hamilton cycle.  When lifted back to $G_{n,p}$ via Theorem \ref{thrio} we get the promised $K_r$-cycle.

If $p=\om^{-1}n^{-2/r}\log^{1/\binom{r}{2}}n$ then Lemma 1.4 of \cite{JKV} implies that w.h.p. there will be vertices that are not in a copy of $K_r$.
\qed

This completes the proof of Theorem \ref{th1}.
\section{Discussion and open problems}
We first note that we can replace $K_r$ by any strictly 1-balanced graph $F$ and then apply Theorem 15 of \cite{Rio} and obtain a spanning subgraph made up of a sequence of edge disjoint copies of $F$, where adjacent copies in the sequence share exactly one common vertex. More precisely, for a graph $F$ we let $d_1(F)=\frac{|E(F)|}{|V(F)|-1}$. A graph is strictly 1-balanced if $d_1(F)>d_1(F')$ for all subgraphs $F'\subseteq F$ with at least two vertices. Theorem 15 amends Theorem \ref{thrio} by having the requirement that $p\leq n^{-1/d_1+\e}$ and letting $\p=ap^{|E(F)|}$ for some constant $a>0$. Note that $|E(F)|=\binom{r}{2}$ $d_1(K_r)=r/2$ and so Theorem \ref{th1} is just a special case, other than the knowledge that we can take $a=1$.
We call the constructions that arise $F$-cycles.

There is a weakness in the result. Consider the diagram below:
\begin{center}
\def\C4{
\draw [fill=black] (0,0) circle [radius=.05];
\draw [fill=black] (0,1) circle [radius=.05];
\draw [fill=black] (0.5,0.5) circle [radius=.05];
\draw [fill=black] (-0.5,0.5) circle [radius=.05];
\draw (0,0) -- (0.5,0.5);
\draw (0,0) -- (-0.5,0.5);
\draw (0,1) -- (0.5,0.5);
\draw (0,1) -- (-0.5,0.5);
}
\pic{
\begin{scope}[shift={(0,0)}, scale=1]
\C4
\end{scope}
\begin{scope}[shift={(1,0)}, scale=1]
\C4
\end{scope}
\begin{scope}[shift={(2,0)}, scale=1]
\C4
\end{scope}
\begin{scope}[shift={(3,0)}, scale=1]
\C4
\end{scope}
\begin{scope}[shift={(3,1)}, scale=1]
\C4
\end{scope}
\begin{scope}[shift={(3,2)}, scale=1]
\C4
\end{scope}
\begin{scope}[shift={(2,2)}, scale=1]
\C4
\end{scope}
\begin{scope}[shift={(1,2)}, scale=1]
\C4
\end{scope}
\begin{scope}[shift={(0,2)}, scale=1]
\C4
\end{scope}
\begin{scope}[shift={(0,1)}, scale=1]
\C4
\end{scope}
\node at (1.5,-1) {$C_4$-cycle};
}
\end{center}
We cannot use the above argument to show that the threshold for an $n$-vertex copy of the above example has a threshold at $p=n^{-3/4+o(1)}$. The reason being that we have no control over the positioning of the connecting vertices i.e. we cannot prevent something like the following being part of the $F$-cycle:
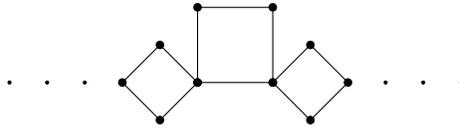
\begin{center}
\def\C4{
\draw [fill=black] (0,0) circle [radius=.05];
\draw [fill=black] (0,1) circle [radius=.05];
\draw [fill=black] (0.5,0.5) circle [radius=.05];
\draw [fill=black] (-0.5,0.5) circle [radius=.05];
\draw (0,0) -- (0.5,0.5);
\draw (0,0) -- (-0.5,0.5);
\draw (0,1) -- (0.5,0.5);
\draw (0,1) -- (-0.5,0.5);
}
\def\K4{
\draw [fill=black] (0,0) circle [radius=.05];
\draw [fill=black] (0,1) circle [radius=.05];
\draw [fill=black] (1,0) circle [radius=.05];
\draw [fill=black] (1,1) circle [radius=.05];
\draw (0,0) -- (0,1);
\draw (0,1) -- (1,1);
\draw (1,1) -- (1,0);
\draw (1,0) -- (-0,0);
}
\pic{
\draw [fill=black] (-2,0.5) circle [radius=.02];
\draw [fill=black] (-1,0.5) circle [radius=.02];
\draw [fill=black] (-1.5,0.5) circle [radius=.02];
\begin{scope}[shift={(0,0)}, scale=1]
\C4
\end{scope}
\begin{scope}[shift={(0.5,0.5)}, scale=1]
\K4
\end{scope}
\begin{scope}[shift={(2,0)}, scale=1]
\C4
\end{scope}
\draw [fill=black] (4,0.5) circle [radius=.02];
\draw [fill=black] (3.5,0.5) circle [radius=.02];
\draw [fill=black] (3,0.5) circle [radius=.02];
}
\end{center}
It is therefore an open question as to the threshold for the existence of a spanning $C_4$-cycle. 

The proof also breaks if our adjacent copies share two or more vertices, as in the diagrams below:
\begin{center}
\def\K4{
\draw [fill=black] (0,0) circle [radius=.05];
\draw [fill=black] (0,1) circle [radius=.05];
\draw [fill=black] (1,0) circle [radius=.05];
\draw [fill=black] (1,1) circle [radius=.05];
\draw (0,0) -- (0,1);
\draw (0,1) -- (1,1);
\draw (1,1) -- (1,0);
\draw (1,0) -- (-0,0);
}
\pic{
\begin{scope}[shift={(0,0)}, scale=1]
\K4
\end{scope}
\begin{scope}[shift={(1,0)}, scale=1]
\K4
\end{scope}
\begin{scope}[shift={(2,0)}, scale=1]
\K4
\end{scope}
\begin{scope}[shift={(3,0)}, scale=1]
\K4
\end{scope}
\begin{scope}[shift={(3,1)}, scale=1]
\K4
\end{scope}
\begin{scope}[shift={(3,2)}, scale=1]
\K4
\end{scope}
\begin{scope}[shift={(2,2)}, scale=1]
\K4
\end{scope}
\begin{scope}[shift={(1,2)}, scale=1]
\K4
\end{scope}
\begin{scope}[shift={(0,2)}, scale=1]
\K4
\end{scope}
\begin{scope}[shift={(0,1)}, scale=1]
\K4
\end{scope}
\node at (2,-1) {$C_4$-cycle, overlap 2};
}
\def\K4-{
\draw [fill=black] (0,0) circle [radius=.05];
\draw [fill=black] (0,1) circle [radius=.05];
\draw [fill=black] (1,0) circle [radius=.05];
\draw [fill=black] (1,1) circle [radius=.05];
\draw (0,0) -- (0,1);
\draw (0,1) -- (1,1);
\draw (1,1) -- (1,0);
\draw (1,0) -- (-0,0);
\draw (0,0) -- (1,1);
}
\hspace{1in}
\pic{
\begin{scope}[shift={(0,0)}, scale=1]
\K4-
\end{scope}
\begin{scope}[shift={(1,0)}, scale=1]
\K4-
\end{scope}
\begin{scope}[shift={(2,0)}, scale=1]
\K4-
\end{scope}
\begin{scope}[shift={(3,0)}, scale=1]
\K4-
\end{scope}
\begin{scope}[shift={(3,1)}, scale=1]
\K4-
\end{scope}
\begin{scope}[shift={(3,2)}, scale=1]
\K4-
\end{scope}
\begin{scope}[shift={(2,2)}, scale=1]
\K4-
\end{scope}
\begin{scope}[shift={(1,2)}, scale=1]
\K4-
\end{scope}
\begin{scope}[shift={(0,2)}, scale=1]
\K4-
\end{scope}
\begin{scope}[shift={(0,1)}, scale=1]
\K4-
\end{scope}
\node at (2,-1) {($K_4-e$)-cycle, overlap 2};
}

\end{center}
One can check that the probability an edge occurs in $H$ is not sufficient to imply the existence of a Hamilton cycle of the requisite type as in \cite{DFLS}. For the first example, the expected number of copies of a spanning $C_4$-cycle in $G_{n,p}$ is given by $n!p^{3n/2}$ and so we should take $p\approx n^{-2/3}$. But then $\p$ will be chosen as $\approx n^{-8/3}$ and this is below the threshold of $\om n^{-2}$ for a Hamilton cycle of the required type, see Theorem 3(iii) of \cite{DF1}. We have a similar experience with the second example, with $p\approx n^{-1/2}$ and $\p\approx n^{-5/2}$.

On the other hand, a recent result of Frankston, Kahn, Narayanan and Park \cite{FKNP} enables us to argue that the suggested thresholds are no worse than $\log n$ from the correct values.

\end{document}